




\input amstex
\documentstyle{amsppt}
\magnification=\magstep1
\pagewidth{6.2truein}
\pageheight{8.3truein}
\NoBlackBoxes
\def\chix{{\raise.5ex\hbox{$\chi$}}}
\def\ep{\varepsilon}
\def\A{{\Cal A}}
\def\E{{\Cal E}}
\def\L{{\Cal L}}
\def\S{{\Cal S}}
\def\esssup{\mathop{\text{\rm ess sup}}\limits}
\def\sgn{\operatorname{sgn}}
\def\comp{{\Bbb C}}
\def\tee{{\Bbb T}}
\def\defeq{\mathop{\ \buildrel \text{df}\over =\ }\nolimits}
\topmatter
\title On an inequality of A.~Grothendieck concerning operators on $L^1$
\endtitle
\author Haskell Rosenthal\endauthor
\affil Department of Mathematics\\ The University of Texas at Austin\\
Austin, TX 78712\endaffil
\email rosenthl@math.utexas.edu\endemail
\abstract 
In 1955, A.~Grothendieck proved a basic inequality which shows that any 
bounded linear operator between $L^1(\mu)$-spaces maps (Lebesgue-) dominated 
sequences to dominated sequences. 
An elementary proof of this inequality is obtained via a new 
decomposition principle for the lattice of measurable functions. 
An exposition is also given of the M.~L\'evy extension theorem for operators 
defined on subspaces of $L^1(\mu)$-spaces.
\endabstract 
\endtopmatter

\document
\baselineskip=18pt

\head 1. Introduction\endhead 

Let $\mu,\nu$ be measures on measurable spaces, and let $T:L^1(\mu)\to 
L^1(\nu$ be a bounded linear operator (here $L^1(\mu)$ denotes the real or 
complex Banach space of (equivalence classes of) $\mu$-integrable functions). 
In \cite{G}, (see Corollaire, page 67)
Grothendieck establishes the following fundamental inequality: 
$$\left\{\eqalign{
&\text{Given $f_1,\ldots,f_n$ in $L^1(\mu)$, then}\cr 
&\int \max_i |Tf_i|\,d\nu \le \|T\| \int \max_i |f_i|\,d\mu\ .\cr}
\right.
\tag 1$$ 

We first give some motivation for the inequality, then give a proof 
involving an apparently new principle concerning the lattice of 
measurable functions. 

It follows easily from (1) that every such operator maps {\it dominated\/} 
(or {\it order bounded\/}) sequences into dominated sequences. 
In fact, it follows that 
$$\left\{ \eqalign{
&\text{\it if $F$ is a family in $L^1(\mu)$ for which there exists a 
$\mu$-integrable $\varphi$ with}\cr
&\text{\it $|f|\le\varphi$ {\rm a.e.} for all $f$ in $F$, then there 
exists a non-negative $\nu$-integrable}\cr
&\text{\it  $\psi$ with $\int \psi\, d\nu\le \|T\| \int \varphi\,d\mu$ 
so that $|Tf|\le \psi$ {\rm a.e.} for all $f$ in $F$.}\cr}
\right.
\tag 2$$
This consequence of (1) (which is of course equivalent to (1)) is drawn 
explicitly by Grothendieck in \cite{G} (see Proposition~10, page 66). 

In the summer of 1979, during her research visit to the University of 
Texas at Austin, I suggested to Mireille L\'evy that the inequality (1) 
might actually characterize those operators from a subspace of $L^1(\mu)$ 
to $L^1(\nu)$, which extend to an operator on all of $L^1(\mu)$. 
She indeed confirmed my conjecture \cite{L}. 
Combining L\'evy's result with (1) and a simple application of the 
closed graph theorem, we obtain the 

\proclaim{Extension Theorem} 
Let $\mu,\nu$ be measures on measurable spaces, $X$ a closed linear 
subspace of $L^1(\mu)$, and $T: X\to L^1(\nu)$ a bounded linear operator. 
Then the following assertions are equivalent: 
\roster
\item"(a)" $T$ maps dominated sequences to dominated sequences.
\item"(b)" There is a constant $C$ so that 
$$\left\{ \eqalign{
&\text{given $n$ and $f_1,\ldots,f_n$ in $X$, then}\cr 
&\int \max_i |Tf_i|\,d\nu \le C\int \max_i |f_i|\,d\mu\ .\cr}
\right. 
\tag 3$$
\item"(c)" There is a bounded linear operator $\tilde T:L^1(\mu)\to 
L^1(\nu)$ with $\tilde T|X=T$.
\endroster
Moreover if $\alpha$ denotes the smallest $C$ satisfying $(3)$, then 
$\tilde T$ may be chosen with $\|\tilde T\|=\alpha$.
\endproclaim 

A remarkable development of the setting for the Extension Theorem 
has recently been given in a series of papers by G.~Pisier. 
In \cite{P1}, Pisier obtains an extension theorem 
for operators on $H^1$ to $L^1(\mu)$ which are also bounded from $H^\infty$ 
to $L^\infty (\mu)$. 
In \cite{P2, Theorem~3}, he obtains the appropriate generalization of the 
Extension Theorem for operators from a subspace of $L^p(\mu)$ to $L^1(\nu)$, 
$1\le p\le\infty$, and in fact in the more general setting of Banach lattices. 
Finally, in \cite{P3, Theorem~3.5}, Pisier obtains a non-commutative version 
of the Theorem.
In Section~3, we give a proof of the Extension Theorem following the 
approach in \cite{P1}. 
This also yields a rather quick alternate ``functional-analytical'' 
proof of (1). 
For a given subspace $X$ of $L^1$, our exposition yields an explicit 
representation for elements of $X(L^\infty)$, the closure of 
$X\otimes L^\infty$ in $L^1(L^\infty)$ (see the Corollary towards the 
end of Section~3), which also suggests an open question regarding 
$X(L^\infty)$ (see the second Remark following the Corollary's statement). 

We note one last motivating connection.
Grothendieck's ``$L^1$-inequality'' (1) 
follows immediately from the classical Banach lattice 
result that every such operator $T$ has an absolute value, or modulus, 
$|T|$, which is a linear operator from $L^1(\mu)$ to $L^1(\nu)$ with 
$\| (T)\| = \|T\|$ and 
$$|Tf| \le |T|\, |f|\ \text{ for all } f\in L^1(\mu) 
\tag 4$$ 
(cf.\ \cite{S}). 
However the existence of $|T|$ may readily be deduced from (1), 
which thus certainly appears more basic and elementary. 

\head 2. A decomposition principle for the lattice of measurable functions 
\endhead 

We first formulate the principle for the case of real scalars. 

\proclaim{Lemma 1} 
Let $f_1,\ldots,f_n$ be real valued measurable functions on a measurable 
space. 
There  exist $k$ (depending only on $n$) and non-negative measurable 
functions $h_1,\ldots,h_k$ so that 
\roster
\item"(i)" $h_1+\cdots + h_k = |f_1| \vee \cdots \vee |f_n|$\ ;
\item"(ii)" for all $i$, there exist $\ep_{ij} \in \{0,1,-1\}$ with 
$f_i = \sum_{j=1}^k \ep_{ij} h_j$. 
\endroster 
\endproclaim 

\remark{Remark} 
We do not need the fact that the $k$ in Lemma~1 depends only on $n$. 
Nevertheless, let $k(n)$ be the optimal choice for $k$. 
What is $k(n)$? 
The order of magnitude of $k(n)$? 
Shortly after circulating the original version of this paper, V.~Mascioni 
completely solved this problem, proving that one may choose $k(n)=2^n$, 
and this is best possible \cite{M}. 
(Our proof below yields only that $k(n) \le e^{1/2} 2^n n!$; also see 
the remark following Lemma~2.) 
\endremark

We first deduce the Grothendieck inequality for real scalars from Lemma~1. 
Given $T$ and $f_1,\ldots,f_n$ in $L^1(\mu)$, choose $h_1,\ldots,h_k$ and 
the $\ep_{ij}$'s as in the Lemma. 
Then for each $i$, we have 
$$|Tf_i| = \Big| \sum \ep_{ij} Th_j\Big| \le \sum_j |Th_j|\ .
\tag 5$$ 
Hence 
$$\max_i |Tf_i| \le \sum_j |Th_j|\ .
\tag 6$$
Thus 
$$\align 
\int \max_i |Tf_i|\,d\nu & = \int \sum_j |Th_j|\,d\nu \qquad
\text{ by (6)}\tag 7\cr 
&= \sum_j \int |Th_j|\,d\nu \cr 
&\le \|T\| \sum_j \int h_j\,d\mu \qquad
\text{ since } h_j \ge0 \text{ for all } j\cr
&= \|T\| \int \sum_j h_j\,d\mu  \cr 
&= \|T\| \int \max |f_i|\,d\mu\qquad  \text{ by (i) of the Lemma.} 
\endalign$$

\demo{Proof of Lemma 1} 

We prove the result by induction on $n$. 
Let $(\Omega, \S)$ be the associated measurable space; i.e., $\S$ is a 
$\sigma$-algebra of subsets of $\Omega$, and the $f_i$'s are $\S$-measurable 
functions defined on $\Omega$. 
For $n=1$, let $h_1 = f_1^+$, $h_2 = f_1^-$ 
(where as usual, e.g., $f_1^+ (\omega) = f_1(\omega)$ if $f_1(\omega)\ge0$; 
$f_1^+ (\omega) =0$ otherwise). 
Of course then $|f_1| = h_1+h_2$, $f_1 = h_1 -h_2$. 
Now let $n\ge1$, and suppose the Lemma proved for $n$. 
Let $f_1,\ldots,f_{n+1}$ be given measurable functions on $\Omega$. 
Choose disjoint measurable sets $E_1,\ldots,E_{n+1}$ so that 
$\Omega = \bigcup_{i=1}^{n+1} E_i$ and 
$$|f_1|(\omega) \vee \cdots \vee |f_{n+1}| (\omega ) 
= |f_i(\omega)| \text{ for all $\omega \in E_i$, all $i$.} 
\tag 8$$
Now fix $i$ and apply the induction hypothesis to $f_1,\ldots, f_{i-1}$, 
$f_{i+1},\ldots,f_{n+1}$ on $E_i$. 
We obtain $h_{i1},\ldots,h_{ik} \ge 0$ ($k$ depends only on $n$) 
measurable functions so that 
$$\sum_{j=1}^k h_{ij} = (|f_1| \vee \cdots \vee |f_{i-1}| \vee 
|f_{i+1}| \vee \cdots \vee |f_{n+1}|) \chix_{E_i} 
\ \buildrel \text{df}\over =\ 
\tau_i 
\tag 9$$ 
and so that 
for each $j\ne i$, there are numbers $\ep_{j\ell}^i$ in $\{0,1,-1\}$ with  
$$f_j \chix_{E_i} = \sum_{\ell=1}^k \ep_{j\ell}^i h_{i\ell} \ .
\tag 10$$
Let $E_i^+ = \{\omega : f_i (\omega) \ge 0\}$. 
$E_i^- = \{\omega :f_i (\omega) <0\}$. 
We now claim the following family of functions works, for our ``$h_i$'s'' 
for $n+1$: 
$$\left\{ \eqalign{ 
&h_{i\ell} \chix_{E_i^+} \ ,\ h_{i\ell} \chix_{E_i^-} \ ,\cr 
&(f_i - \tau_i) \chix_{E_i^+}\ ,\ (-f_i-\tau_i) \chix_{E_i^-} \cr 
&(1\le i\le n+1\ ,\ 1\le \ell\le k)\ .\cr} \right. 
\tag 11$$ 
Evidently if $k'$ denotes the total number of functions listed in (11), then 
$$k'= 2(n+1)(k+1)\ .
\tag 12$$
Now, all of these functions are non-negative (the last two types because 
$|f_i| \ge \tau_i$ on $E_i$, by (8)). 
To verify (i) of the Lemma, note that for each $i$, 
$$|f_i| \chix_{E_i} = (f_i-\tau_i)\chix_{E_i^+} + \sum_{\ell=1}^k h_{i\ell} 
\chix_{E_i^+} + (-f_i-\tau_i)\chix_{E_i^-} + \sum_{\ell=1}^k h_{i\ell} 
\chix_{E_i^-}\ .
\tag 13$$ 
Thus, letting  $h_1,\ldots,h_{k'}$ be the functions listed 
in (11), we have that 
$$|f_1| \vee \cdots \vee |f_{n+1}| = \sum_{i=1}^{n+1} |f_i|\chix_{E_i} 
= \sum_{r=1}^{k'} h_r\ .
\tag 14$$ 
Finally, to verify (ii), fix $j$. 
Then 
$$\leqalignno{
f_j \chix_{E_j} & = f_j \chix_{E_j^+} + f_j \chix_{E_j^-} &(15)\cr 
& = (f_j-\tau_j)\chix_{E_j^+} + \tau_j \chix_{E_j^+} - (-f_j-\tau_j) 
\chix_{E_j^-} - \tau_j \chix_{E_j^-} \cr 
& = (f_j -\tau_j) \chix_{E_j^+} + \sum_{\ell=1}^k h_{j\ell} \chix_{E_j^+} 
- (-f_j-\tau_j) \chix_{E_j^-} + \sum_{\ell=1}^k - h_{j\ell} \chix_{E_j^-}\ .
\cr}$$ 
Thus from (10) and (15), we obtain 
$\ep_{jr} = 0,1$, or $-1$ for all $r$ so that 
$$\displaylines{(16)\hfill
f_j = \sum_{i=1}^{n+1} f_j \chix_{E_i} = \sum_{r=1}^{k'} \ep_{jr} h_r\ .
\hfill\qed}$$
\enddemo

We next treat the case of complex scalars. 

\proclaim{Lemma 2} 
Let $f_1,\ldots,f_n$ be complex valued measurable functions on a 
measurable space. 
There exist $k$ (depending only on $n$) and non-negative measurable functions 
$h_1,\ldots,h_k$ so that 
\roster 
\item"(i)" $h_1+\cdots + h_k = |f_1| \vee \cdots \vee |f_n|$. 
\item"(ii)" for all $i$, there exist measurable functions $\ep_{ij}$ 
so that $|\ep_{ij}(\omega)| =1$ or $0$ for all $j$, with 
$f_i = \sum_{j=1}^k \ep_{ij} h_j$
\endroster
\endproclaim 

\remark{Remark}
Let $k_\comp (n)$ denote the optimal choice for $k$.
As in the real scalars case, we again ask what is the order of 
magnitude of $k_\comp$? Our argument below yields that $k_\comp (n)\le en!$.
(V.~Mascioni has also solved this problem, proving that $k_\comp (n)=2^n-1$ 
\cite{M}.) 
\endremark 

The deduction of the Grothendieck $L^1$-inequality involves the following 

\proclaim{Corollary} 
Let $f_1,\ldots,f_n$ be as in Lemma 2, and let $\ep>0$.  
There exist $h_1,\ldots, h_k$ non-negative measurable functions satisfying 
(i) of Lemma~1 and 

\noindent 
(ii) for all $i$ there exist numbers $\alpha_{ij}$ with $|\alpha_{ij}|=1$ 
or $0$ for all $j$, and  
$$\Big|f_i - \sum_{j=1}^k \alpha_{ij} h_j\Big| \le \ep 
(|f_1| \vee \cdots \vee |f_n|)\ .
\tag 17$$
\endproclaim 

\remark{Comment} 
If the $\ep_{ij}$'s in Lemma 2 can be chosen as simple functions (which is of 
course the case if the $f_i$'s are simple), then the dependence of the 
$\alpha_{ij}$'s on $\ep$ may be eliminated; i.e., we then have 
$f_i = \sum_j \alpha_{ij} h_j$ for all $i$. 
Note this is the case if the $f_i$'s are all real-valued; thus Lemma~2 
implies Lemma~1. 
\endremark 

\demo{Proof of the Corollary using Lemma 2} 
Let the $h_i$'s and $\ep_{ij}$'s satisfy the conclusion of Lemma~2. 
We may choose disjoint measurable sets $F_1,\ldots,F_r$ with 
$\Omega = \bigcup_{i=1}^r F_i$, so that for every $\nu$, 
$1\le \nu \le r$, every $i$, $1\le i\le n$, and all $j$, $1\le j\le k$, 
there is a number $\ep_{ij}^\nu$, with $|\ep_{ij}^\nu | =1$ or 
$\ep_{ij}^\nu =0$, so that 
$$|\ep_{ij}(\omega) - \ep_{ij}^\nu| \le \ep \text{ for all }\omega \in F_\nu\ .
\tag 18$$ 
We now claim: {\it The family of functions 
$$h_i\chix_{F_\nu}\qquad 1\le i\le k\ ,\ 1\le \nu \le r\ ,$$ 
serves as our ``$h_\ell$'s''; for each $i$, the constant $\ep_{ij}^\nu$ 
serves as our ``$\alpha_{i\ell}$.''}
Indeed, we have that 
$$\sum_{i,\nu} h_i \chix_{F_\nu} = |f_1| \vee \cdots \vee |f_n|\ .
\tag 19$$ 
Finally, fix $i,\nu$. 
Then 
$$\leqalignno{
\Big|f_i \chix_{F_\nu} - \sum \ep_{ij}^\nu h_j\chix_{F_\nu}\Big| 
& = \Big| \sum_j (\ep_{ij} -\ep_{ij}^\nu) h_j \chix_{F_\nu}\Big| 
\ \text{ by Lemma 2(ii)}&(20)\cr
&\le \ep \sum h_j \chix_{F_\nu} 
= \ep (|f_1| \vee \cdots \vee |f_n|) \chix_{F_\nu}\ .\cr}$$ 
Since the $F_\nu$'s are a partition of $\Omega$, the result is proved.\qed
\enddemo 

\demo{Proof of Lemma 2} 
Again we proceed by induction. 
For any measurable complex valued function $f$, let 
$$(\sgn f) (\omega) = {f(\omega)\over |f(\omega)|} \ \text{ if } 
f(\omega) \ne0\ ,\ \
\sgn f(\omega) =0 \text{ otherwise.}$$ 
Of course now the $n=1$ case is ``completely'' trivial; simply let 
$h_1=|f_1|$ and 
$$\ep_1(\omega) = \sgn f_1(\omega)\ .$$
Again, suppose Lemma 2 proved for $n$, and let $f_1,\ldots,f_{n+1}$ be 
given measurable functions. 
Choose the measurable partition $E_1,\ldots,E_{n+1}$ satisfying (8), 
and proceed exactly as in the case of Lemma~1. 
Thus, we obtain $h_{ij}$'s, $1\le j\le k$ satisfying (9) (with $\tau_i$ 
as defined in (9)), so that for each $j\ne i$, there are measurable functions 
$\ep_{j\ell}^i$ with $|\ep_{j\ell}^i (\omega)| =0$ or $1$ for all $\omega$, 
satisfying (10). 
Now we claim that the family of ``$h_i$'s'' may be taken to be 
$$h_{i\ell} \ ,\quad (|f_i|-\tau_i)\chix_{E_i}\ ,\quad
1\le i\le n+1\ ,\ 1\le \ell\le k\ .
\tag 21$$ 
Thus listing these as $h_1,\ldots, h_{k'}$, we  have 
$$k'= (k+1)(n+1) \ .
\tag 22$$ 
Lemma 2(i) now follows immediately, for 
$$|f_i| \chix_{E_i} = (|f_i| -\tau_i) \chix_{E_i} + \sum_{\ell=1}^k h_{i\ell} 
\ \text{ for all }\ i\ .
\tag 23$$ 

It remains to verify (ii). Fix $j$. 
Then 
$$(f_j-(\sgn f_j) \tau_j) \chix_{E_j} 
= (\sgn f_j) (|f_j| -\tau_j)\chix_{E_j}\ . 
\tag 24$$ 
Thus 
$$f_j \chix_{E_j} = (\sgn f_j)(|f_j|-\tau_j\chix_{E_j}) 
+\sum_{\ell=1}^k (\sgn f_j) h_{j\ell} \quad \text{ by (9).}
\tag 25$$
Combining (10) and (25), we thus obtain our measurable functions 
$\ep_{j1},\ldots,\ep_{jk'}$ valued in $\tee\cup\{0\}$ with 
$$\leqalignno{f_j & = \sum_{i\ne j} f_j \chix_{E_i} + f_j\chix_{E_j}&(26)\cr 
&= (\sgn f_j) (|f_j| -\tau_j \chix_{E_j}) + \sum_{\ell=1}^k (\sgn f_j) h_{j
\ell} + \sum_{i\ne j} \ \sum_{\ell=1}^k \ep_{j\ell}^i h_{i\ell}\cr 
& = \sum_{\ell=1}^{k'} \ep_{j\ell} h_\ell\ .\cr}$$
\rightline{$\square$}
\enddemo 

We conclude Section 2 
with a deduction of the complex Grothendieck $L^1$-inequality. 
Let then $\mu,\nu$ be measures on measurable spaces, $T: L^1 (\mu)\to L^1(\nu)$ 
be a bounded linear operator, and $f_1,\ldots,f_n$ in $L^1(\mu)$ be given. 
Let $\ep>0$ be given, and choose $h_1,\ldots,h_k$ and the complex numbers 
$\alpha_{ij}$ as in the conclusion of the Corollary to Lemma~2. 

Now for each $i$, define $p_i$ by 
$$p_i = f_i - \sum_{j=1}^k \alpha_{ij} h_j\ .
\tag 27$$ 
Then we have that $f_i = \sum \alpha_{ij} h_j + p_i$  and moreover 
$$|p_i| \le \ep (|f_1| \vee \cdots \vee |f_n|)\ \text{ by (ii) of the 
Corollary.}
\tag 28$$ 
Thus 
$$\leqalignno{|Tf_i| & = \Big| \sum_j \alpha_{ij} Th_j +Tp_i\Big| &(29)\cr 
&\le \sum_j |Th_j| + |Tp_i| \ \text{ since $|\alpha_{ij}|\le 1$ for all $j$}\cr
&\le \sum_j |Th_j| + \sum_j |Tp_j|\ .\cr}$$
Thus also 
$$\max_i |Tf_i| \le \sum_j (|Th_j | + |Tp_j|)\ ,
\tag 30$$ 
whence
$$\leqalignno{\qquad \int \max_i |Tf_i|\,d\nu
&\le \sum_j \int (|Th_j| + |Tp_j|)\,d\nu &(31)\cr 
&\le \|T\| \biggl( \int \sum_j h_j \,d\mu + \int \sum |p_j|\,d\mu\biggr)\cr 
&\le (1+n\ep)\|T\| \int \max_i |f_i|\,d\mu 
\text{ by (i) of the Corollary and (28).}
\cr}$$
Since $\ep>0$ is arbitrary, the inequality (1) is proved.\qed

\head 3. A proof of the Extension Theorem\endhead 

As noted in the introduction, we follow the approach in \cite{P1}, thus 
obtaining an alternate proof of the Grothendieck $L^1$-inequality. 
(The approach, despite its brevity, seems considerably more sophisticated 
than the elementary proof given by our decomposition result, however.) 
Throughout, let $\mu,\nu$ and $T$ be as in the statement of the 
Extension Theorem. 
We shall also assume that $\nu$ is ``nice enough'' so that  
$(L^1(\nu))^* = L^\infty (\nu)$ (any $L^1(\nu)$ is isometric to 
$L^1(\nu')$ with $\nu'$ nice). 

(a) $\Rightarrow$ (b) 
For $Y$ a subspace of $L^1(\mu)$ or 
$L^1(\nu)$, let $Y_d$ denote the space of all dominated sequences $(y_n)$ 
in $Y$, under the norm $\| (y_n)\|_d = \int\sup_n |y_n|\,d\mu$. 
We easily check that $Y_d$ is a Banach space; evidently then $T$ induces 
a linear operator $S$ from $Y_d$ to $(L^1(\nu))_d$, which has closed graph, 
since $T$ itself is bounded. 
Thus $S$ is bounded. 

(c) $\Rightarrow$ (a) follows immediately from Grothendieck's 
$L^1$-inequality (1). 
We give here an alternate proof of (1), using the set up in \cite{P1}. 
We will freely use here some standard facts about $Y\widehat\otimes Z$, 
the projective tensor product of Banach spaces $Y$ and $Z$. 
Let $L^1 (\mu,Y)$ denote the space of Bochner-integrable $Y$-valued 
functions on $\Omega$ (where $(\Omega,\S,\mu)$ is the measure space 
associated to $\mu$). 
{\it Then $L^1(\mu,Y)$ is $($canonically isometric to\/})  
$L^1(\mu) \widehat\otimes Y$ (see Th\'eor\`eme~2, page~59 of \cite{G}). 
It follows immediately that $T\otimes I$ yields a linear operator from 
$L^1(\mu)\otimes Y$ to $L^1(\nu)\otimes Y$ with $\|T\otimes I\|=\|T\|$. 
(Here, we assume ``$X$''~$= L^1(\mu)$; i.e., the hypotheses of (1).) 
We apply this fact to $Y=L^\infty (\nu)$. 
It follows that for any $n$, $f_1,\ldots,f_n$ in $L^1(\mu)$, and 
$\varphi_1,\ldots,\varphi_n$ in $L^\infty (\nu)$. 
$$\Big\|\sum_{i=1}^n Tf_i\otimes \varphi_i\Big\| 
\le \|T\|\, \Big\|\sum_{i=1}^n f_i \otimes\varphi_i\Big\|\ .
\tag 32$$ 
Here, $g \defeq \sum f_i\otimes \varphi_i$ denotes the element of 
$L^1(\mu,L^\infty (\nu))$ defined by $g(\omega) = \sum f_i(w) \varphi_i$, 
$\omega\in\Omega$; note that 
$$\|g\| = \int \|g(\omega) \|\,d\mu (\omega) 
= \int \esssup_s \Big| \sum f_i (\omega) \varphi_i (s)\Big|\,d\mu (\omega)\ .
\tag 33$$ 
Now fixing  $f_1,\ldots,f_n$ and $\varphi_1,\ldots,\varphi_n$ as above, 
we have 
$$\leqalignno{
&\Big| \int \sum (Tf_j)(s) \varphi_j(s)\, d\nu(s)\Big| &(34)\cr
&\qquad 
\le \int \esssup_t \Big| \sum_{j=1}^n (Tf_j)(s)\varphi_j(t)\Big|\,d\nu(s)\cr
&\qquad 
\le \|T\| \int \esssup_t \Big| \sum f_j(\omega)\varphi_j(t)\Big|\,d\mu 
(\omega)\qquad \text{ by (32) and (33)}\cr 
&\qquad 
\le \|T\| \biggl( \int \max_j |f_j(\omega)|\,d\mu(\omega)\biggr) 
\Big\| \sum |\varphi_j|\,\Big\|_{L^\infty (\nu)}\ .\cr}$$ 
Now since $\nu$ is nice, a standard argument yields that we may choose 
$\varphi_1,\ldots,\varphi_n$ in $L^\infty (\nu)$ with 
$\|\sum |\varphi_j|\, \|_{L^\infty (\nu)} =1$ and 
$$\int \max |Tf_j|(s)\,d\nu (s) 
= \int \sum (Tf_j)(s) \varphi_j (s)\, d\nu (s)\ .
\tag 35$$
Evidently (34) and (35) immediately yield Grothendieck's inequality~(1). 

It remains to prove (c) $\Rightarrow$ (d) and the ``moreover'' statement, 
i.e., M.~L\'evy's theorem. 
We closely follow the brief sketch given by Pisier in \cite{P1}, 
crystallizing some elements of the discussion. 
It is convenient to introduce one more condition in the Extension Theorem, 
which is explicitly used in \cite{P1}. 
\roster
\item"(d)" {\it There is a constant $C$ so that for any  $n$, 
$f_1,\ldots,f_n$ in $X$, and simple 
$\varphi_1,\ldots,\varphi_n$ in $L^\infty (\nu)$,}
$$\Big|\sum_i  \int (Tf_i)\varphi_i\,d\nu\Big| 
\le C\int \esssup_s \Big| \sum f_i (\omega)\varphi_i(s)\Big|\,d\mu(\omega)\ . 
\tag 36$$
\endroster 

We first prove (d) $\Rightarrow$ (c). 
Consider the following general problem: 
Given Banach spaces $Y$, $B$, $X$ a closed linear subspace of $Y$, 
$T:X\to B^*$ a bounded linear operator, and $C>0$, when does there exist 
$\tilde T:Y\to B^*$ extending $T$, with $\|\tilde T\|\le C$? 
Is there a way of formulating this problem in terms of the Hahn-Banach 
Theorem? 
As e.g., developed in \cite{G}, $\L(Y,B^*)$ is indeed, naturally isometric 
to $(Y\hat\otimes B)^*$. 
The pairing is as follows: given $T:Y\to B^*$ a bounded linear operator 
and $\omega \defeq \sum y_i\otimes b_i$ in $(Y\hat\otimes B)^*$ (with 
$\sum \|y_i\|\, \|b_i\| <\infty$), set 
$$\langle T,\omega\rangle = \sum_i \langle Ty_i,b_i\rangle\ .
\tag 37$$ 
We then obtain the following result: 

\proclaim{Lemma 3}
Given $Y$, $B$, $X$, and $T$ as above, the following are equivalent: 
\roster 
\item"(i)" There is a linear operator $\tilde T:Y\to B^*$ extending $T$, with 
$\|\tilde T\| \le C$.  
\item"(ii)" Let $X_0$, $B_0$ be dense linear subspaces of $X$ and $B$ 
respectively and regard $X_0\otimes B_0$ as a linear subspace of 
$Y\hat\otimes B$. 
Define $F_T$ on $X_0\otimes B_0$ by $F_T (\omega) = \langle T,\omega\rangle$
for all $\omega$ in $X\otimes B$. 
Then 
$$\|F_T\| \le C\ .
\tag 38$$
\endroster
\endproclaim 

To see this, note that (i) $\Rightarrow$ (ii) is immediate. 
If (ii) holds, let $\tilde F_T$ be a Hahn-Banach extension of $F_T$ to 
$Y\hat\otimes B^*$. 
Now simply let $\tilde T$ be the unique element of $\L(X,B^*)$ satisfying 
$$\langle \tilde T,\omega\rangle = \tilde F_T(\omega) \text{ for some } 
\omega \in Y\hat\otimes B^*\ .
\tag 39$$ 
To obtain (d) $\Rightarrow$ (c) of the Extension Theorem, let $X=X_0$, 
$Y= L^1(\mu)$, $B= L^\infty (\nu)$, and $B_0$ the subspace of 
$B$ consisting of simple functions. 
Now condition (d) simply means that $\|F_T\| \le C$, where $F_T$ is as in 
Lemma~3(ii). 
Thus by Lemma~3, we obtain a linear operator 
$\tilde T:L^1(\mu)\to L^1(\nu)^{**}$ extending $T$ (where of course we 
regard $L^1(\nu) \subset L^1(\nu)^{**}$). 
The proof is completed by observing that there exists a norm-one linear 
projection $P$ from $L^1(\nu)^{**}$ onto $L^1(\nu)$; then $P\circ \tilde T$ 
yields the desired operator extending $T$. 

It remains to show that (b) $\Rightarrow$ (d). 

The argument for this implication involves a critical identification, 
due to M.~L\'evy \cite{L}, and appears to have been omitted from the sketch 
given in \cite{P1}. 

\proclaim{Lemma 4} 
Let $B_0$ denote the subspace of $L^\infty (\nu)$ consisting of simple 
functions, and let $g\in X \otimes B_0$. 
Then 
$$\|g\| = \min \biggl\{ \int \max_j |f_j|\,d\mu \Big\| \sum_i |\varphi_i|
\Big\|_\infty \biggr\} 
\tag 40$$ 
the minimum taken over all $n$, $f_1,\ldots,f_n$ in $X$, and 
$\varphi_1,\ldots,\varphi_n$ in $B_0$ so that $g= \sum f_j\otimes\varphi_j$ 
(where  $\|g\|$ is defined as in $(33)$).
\endproclaim 

\demo{Proof of Lemma 4} 
Suppose first $g= \sum f_j\otimes \varphi_j$ where  $f_1,\ldots,f_n$ are in 
$L^1(\mu)$, $\varphi_1,\ldots,\varphi_n$ are in $L^\infty (\nu)$ (we do not 
need to assume here that the $f_i$'s belong to $X$). 
We then have that for any $\omega$ and any $s$, 
$$\Big| \sum_{j=1}^n f_j (\omega)\varphi_j (s)\Big| 
\le \max_j |f_j(\omega)| \sum_j |\varphi_j|(s)\ .
\tag 41$$ 
It follows immediately that 
$$\|g\| \le \int \max |f_j(\omega)|\,d\mu (\omega) \Big\|\sum |\varphi_j|
\Big\|_\infty\ .
\tag 42$$ 
Thus
$$\align
\|g\| &\le \inf \biggl\{ \int \max_j |f_j|\,d\mu \Big\| \sum |\varphi_j|
\Big\|_\infty : g = \sum_{j=1}^n f_j\otimes \varphi_j\tag 43\\ 
&\qquad
\text{with }f_j\in L^1(\mu)\text{ and }\varphi_j \in B_0 \text{ for all }j
\biggr\}\ .
\endalign$$ 
Now $g= \sum_{i=1}^\ell x_i\otimes\psi_i$ with the $x_i$'s in $X$ and the 
$\psi_i$'s in $B_0$. 
We may then choose a $\nu$-measurable partition 
$E_1,\ldots,E_m$ of $S$ so that the $\psi_i$'s are all $\A$-measurable, 
where $\A$ is the algebra generated by the disjoint sets 
$E_1,\ldots,E_n$. 
(Here, we assume $L^1(\nu) = L^1 (S,\E,\nu)$.) 
It {\it then\/} follows that we may choose $z_1,\ldots,z_m$ in $X$ with 
$$g= \sum_{i=1}^m z_i \otimes \chix_{E_i}\ . 
\tag 44$$ 
But then if $\omega\in \Omega$ and $s\in E_i$, 
$$|g(\omega)(s)| = |z_i(\omega)|\ .
\tag 45$$ 
This shows 
$$\|g(\omega)\|_{L^\infty (\mu)} = \max_i |z_i(\omega)|\ .$$ 
Hence 
$$\align \|g\| &= \int \max_i |z_i|\,d\mu\tag 46\\ 
&= \int \max_i |z_i|\,d\mu \Big\| \sum |\chix_{E_j}|\Big\|_\infty \ ,
\endalign$$ 
proving (40).\qed

We finally show that (b) $\Rightarrow$ (d), thus completing the proof 
of the Extension Theorem. 
(The moreover assertion follows from the proof that (d) $\Rightarrow$ (c), 
for of course we show the same constant $C$ in (b) works for (d).) 

Let then $f_1,\ldots,f_n$ be given in $X$, $\varphi_1,\ldots,\varphi_n$ 
be simple elements of $L^\infty (\nu)$, and let $C$ be as in (b). 

By Lemma 4, we may choose $x_1,\ldots,x_m$ in $X$ and $\psi_1,\ldots,\psi_m$ 
simple in $L^\infty (\nu)$ so that letting $g=\sum f_i\otimes \varphi_i$, 
then 
$$\leqalignno{
&g = \sum x_i \otimes \psi_i &\text{(46)(i)}\cr 
&\|g\| = \int \max_i |x_i|\, d\mu \Big\|\sum |\psi_j|\Big\|_\infty\ .
&\text{(46)(ii)}\cr}$$
Now 
$$\align 
\Big| \sum_i \int (Tf_i)\varphi_i\,d\mu\Big| 
& = \Big| \sum_i \langle Tf_i,\varphi_i\rangle \Big| 
= \Big| \sum_i \langle Tx_i,\psi_i\rangle\Big|\text{ (by (46)(i))}\tag 47\cr 
&\le \int \sum_i |Tx_i(s)\psi_i (s)|\, d\nu (s) \cr
&\le \int \max_i |Tx_i|(s) \Big\|\sum_j |\psi_j|\Big\|_\infty\,d\nu(s)\cr
&\le C\int \max |x_i(\omega)|\,d\mu (\omega) 
\Big\| \sum |\psi_j|\Big\|_\infty \text{ (by (b))}\cr 
&= C\|g\|\cr 
&= C\int \esssup_s \Big| \sum f_i(\omega)\varphi_i(s) \Big|\,d\mu (\omega)\ .
\endalign$$ 
This completes the proof of the Extension Theorem.\qed
\enddemo 

The following representation result follows from the above proof of 
M.~L\'evy's theorem, and seems to be what's ``really going on'' (see 
also Lemma~1 of \cite{L}).

\proclaim{Corollary} 
Let $X$ be a closed linear subspace of $L^1(\mu)$, and let 
$X(L^\infty (\nu))$ denote the closure of $X\otimes L^\infty (\nu)$ in 
$L^1(\mu,L^\infty (\nu))$. 
Then given $g\in X(L^\infty (\nu))$ and $\ep>0$, there exists a dominated 
sequence $(x_j)$ in $X$ and a sequence $(\varphi_j)$ in $L^\infty (\nu)$ so 
that 
\roster
\item"(i)" $\sum \varphi_j$ converges unconditionally in $L^\infty (\nu)$.  
\item"(ii)" $g= \sum x_j \otimes \varphi_j$.
\item"(iii)" $\int \sup_j |x_j(\omega)|\,d\mu (\omega) 
\big\|\sum_i |\varphi_i|\big\|_{L^\infty (\nu)} \le \|g\| +\ep$. 
\endroster
\endproclaim 

\remark{Remarks} 
1. If $(x_j)$ in $X$ is dominated and $\sum \varphi_j$ in $L^\infty (\nu)$ 
converges unconditionally, 
then $\sum x_j \otimes \varphi_j$ converges unconditionally 
in $L^1 (\mu,L^\infty (\nu))$, to an element of $X(L^\infty(\nu))$. 
Indeed, for any choice of scalars $(\alpha_j)$ with $|\alpha_j|\le 1$ 
for all $j$ and any $k\le \ell$, we have that 
$$\Big\| \sum_{j=k}^\ell \alpha_j x_j\otimes \varphi_j\Big\| 
\le \int \max_j |x_j (\omega)|\,d\mu (\omega) 
\Big\| \sum_{j=k}^\ell |\varphi_j|\Big\|_\infty\ . 
\tag 48$$ 
But $\sum \varphi_j$ converges unconditionally iff 
$$\Big\| \sum_k^\ell |\varphi_j|\Big\|_\infty \to 0\ \text{ as }\ 
k\to \infty \text{ with } \ell \ge k\ .$$ 
Hence $\sum \alpha_j x_j\otimes \varphi_j$ converges by (48). 

2. Suppose $(x_j)$ in $X$ and $(\varphi_j)$ in $L^\infty (\nu)$ 
satisfy 
$$\int \sup_j |x_j|\,d\mu \Big\|\sum_i |\varphi_i|\Big\|_\infty 
\defeq \tau <\infty\ .$$ 
(Equivalently, $(x_j)$ is dominated and $\sum \varphi_j$ is weakly 
unconditionally summing in $L^\infty (\nu)$.) 
It then follows that for $\mu$-almost all $\omega$, $\sup_j|x_j (\omega)| 
<\infty$; for each such $\omega$, we obtain that $\sum x_j(\omega)\varphi_j$ 
converges absolutely pointwise a.e.\ to an element of $L^\infty (\nu)$, 
and the function $g(\omega) \defeq \sum x_j (\omega)\varphi_j$ belongs 
to $L^1 (\mu,L^\infty (\nu))$ with $\|g\| \le \tau$. 
Does it then follow that $g$ belongs to $X(L^\infty (\nu))$?
This is indeed so provided $X$ is isomorphic to a separable dual space, 
or more generally, a dual space with the Radon-Nikodym property. 
\endremark 

\demo{Proof of the Corollary} 
Letting $B_0$ denote the space of the simple $L^\infty (\nu)$ functions as 
above, we have that $L^1(\mu)\otimes B_0$ is dense in $L^1(\mu)\hat\otimes 
L^\infty (\nu)$ since $B_0$ is dense in $L^\infty (\nu)$. 
Hence given $\ep>0$, we may choose a sequence $(g_j)$ in $L^1(\mu)\otimes 
B_0$ with 
$$\biggl( \sum \|g_j\|^{1/2}\biggr)^2 < \|g\| +\ep\ \text{ and }\  
g=\sum g_j\ . 
\tag 49$$ 
Now by Lemma 4, for each $i$, we may choose finite sequences 
$(x_{ij})_{j=1}^{m_i}$ in $X$ and $(\varphi_{ij})_{j=1}^{m_i}$ in $B_0$ 
with $g_i = \sum_j x_{ij} \otimes \varphi_{ij}$ and 
$$\int \max_j |x_{ij}|\,d\mu (\omega) = \|g_i\|^{1/2} 
= \Big\| \sum_j |\varphi_{ij} |\Big\|_\infty\ .
\tag 50$$ 
Hence the series $\sum_i \sum_{j=1}^{m_i} x_{ij}\otimes \varphi_{ij}$ 
converges unconditionally to $g$. 
Now we have moreover that 
$$\align 
\int \sup_i \max_j |x_{ij}|\,d\mu (\omega) 
&\le \sum_i \int \max_j |x_{ij}|\,d\mu (\omega) \tag 51\cr
&\le \sum |g_j|^{1/2}\ \text{ by (50).}  
\endalign$$ 
Thus the sequence $(x_{ij})$ with $1\le j\le m_i$, $i=1,2,\ldots$ is indeed 
dominated. 
Also $\sum_i \sum_{j=1}^{m_i} \varphi_{ij}$ converges unconditionally in 
$L^\infty (\nu)$ and 
$$\Big\| \sum_i \sum_j |\varphi_{ij}|\Big\|_\infty 
\le \sum_i \Big\|\sum_j |\varphi_{ij}|\Big\|_\infty 
\le \sum |g_j|^{1/2}\text{ by (50).} 
\tag 52$$
The Corollary now follows immediately from (49)--(52).\qed
\enddemo

\enddocument